\documentclass[a4paper,12pt]{amsart}

\usepackage{rotating,pdflscape,color,amssymb,hyperref,subfigure,psfrag,amsmath,eufrak,bbm,epsfig}
\usepackage{a4wide,amscd}
\usepackage{tikz} 
 \usepackage[usenames,dvipsnames]{pstricks}
 \usepackage{pst-grad} 
 \usepackage{pst-plot} 
\usepackage[small]{caption}

\newcommand{\ga}{\gamma}

\newcommand{\ld}{\ldots}
\newcommand{\beg}{\begin}
\newcommand{\en}{\end}

\newcommand{\trm}{\textrm}

\newcommand{\bgt}{\begin{itemize}}
\newcommand{\ent}{\end{itemize}}
\newcommand{\ite}{\item}

\newcommand{\op}{\operatorname}
\newcommand{\eqre}{\eqref}
\newcommand{\re}{\ref}
\newcommand{\la}{\label}

\newcommand{\rfl}{\rfloor}
\newcommand{\lfl}{\lfloor}

\newcommand{\lan}{\langle}
\newcommand{\ran}{\rangle}

\newcommand{\ds}{\displaystyle}

\newcommand{\p}{\mathbb{P}}

\newcommand{\Tr}{\operatorname{Tr}}

\newcommand{\ssi}{if and only if }

\newcommand{\E}{\mathbb{E}}

\newcommand{\R}{\mathbb{R}}
\newcommand{\C}{\mathbb{C}}

\newcommand{\ud}{\mathrm{d}}

\newcommand{\pro}{probability }

\newcommand{\f}{\frac}
\newcommand{\ff}{\frac{1}}
\newcommand{\lf}{\left}
\newcommand{\ri}{\right}

\newcommand{\st}{such that }

\newcommand{\lam}{\lambda}

\newcommand{\ti}{\times}

\newcommand{\vfi}{\varphi}
\newcommand{\ste}{\, ;\, }

\newcommand{\eps}{\varepsilon}

\newcommand{\al}{\alpha}
\newcommand{\tta}{\theta}

\newcommand{\ovl}{\overline}

\newcommand{\bbm}{\begin{bmatrix}}
\newcommand{\ebm}{\end{bmatrix}}
\newcommand{\bes}{\begin{equation*}}
\newcommand{\ees}{\end{equation*}}
\newcommand{\be}{\begin{equation}}
\newcommand{\ee}{\end{equation}}
\newcommand{\beqy}{\begin{eqnarray}}
\newcommand{\eeqy}{\end{eqnarray}}
\newcommand{\beq}{\begin{eqnarray*}}
\newcommand{\eeq}{\end{eqnarray*}}
\newcommand{\one}{\mathbbm{1}}
\newcommand{\lto}{\longrightarrow}

\newcommand{\ie}{\emph{i.e. }}

\newcommand{\bpm}{\begin{pmatrix}}
\newcommand{\epm}{\end{pmatrix}}

\newcommand{\cd}{\cdots}

\newcommand{\bpr}{\beg{proof}}
\newcommand{\epr}{\en{proof}}

\newcommand{\bet}{\beta}
\newcommand{\del}{\delta}

\newcommand{\vb}{\mathbf{v}}

\newcommand{\ha}{\hat{a}}
\newcommand{\ka}{\kappa}
\newcommand{\sL}{\stackrel{sl.}{\sim}}
\newcommand{\bi}{\mathbf{i}}

\newcommand{\bv}{\mathbf{v}}
\newcommand{\br}{\mathbf{r}}

\newcommand{\bfe}{\mathbf{e}}

\newcommand{\tL}{\widetilde{\ell}}

\newtheorem{Th}{Theorem}[section]

\newtheorem{cor}[Th]{Corollary}
\newtheorem{Def}[Th]{Definition}

\date{\today}
\thanks{This work was    
partly accomplished during the first named author's stay at   New York University Abu Dhabi, Abu Dhabi (U.A.E.).}

\title[Heavy tailed band matrices]{Localization and delocalization  for heavy tailed band matrices}
\author{Florent Benaych-Georges} \address{MAP 5, UMR CNRS 8145 - Universit\'e Paris Descartes, 45 rue des Saints-P\`eres 75270 Paris cedex~6, France.} \email{florent.benaych-georges@parisdescartes.fr}

\author{Sandrine P\'ech\'e} \address{LPMA, Universit\'e Paris Diderot, Site Chevaleret 75205 Paris Cedex 13, France.} 
\email{sandrine.peche@math.univ-paris-diderot.fr}

\begin{document}
\maketitle

\beg{abstract}We consider some random band matrices with band-width $N^\mu$ whose entries are independent random variables with distribution tail in $x^{-\al}$. We consider the largest eigenvalues and the associated eigenvectors and prove the following  phase transition. On the one hand, when  $\al<2(1+\mu^{-1})$, the largest eigenvalues have order $N^{\f{1+\mu}{\al}}$, are asymptotically distributed as a Poisson process and their associated eigenvectors are essentially carried  by two coordinates (this phenomenon has already been remarked for full matrices by Soshnikov in \cite{sosh04,sosh06}    when $\al<2
$ and by Auffinger \emph{et al} in \cite{ABP09} when $\al<4$). On the other hand, when  $\al>2(1+\mu^{-1})$, the largest eigenvalues have order $N^{\f{\mu}{2}}$ and most eigenvectors of the matrix are delocalized, \ie approximately uniformly distributed on their $N$ coordinates.
\en{abstract}

\beg{abstract}
On consid\`ere des matrices al\'eatoires \`a structure bande dont la bande a pour largeur $N^\mu$ et dont les coefficients  sont ind\'ependants \`a queue de distribution en $x^{-\al}$. On s'int\'eresse aux plus grandes valeurs propres et aux vecteurs propres associ\'es et prouve la transition de phase suivante. D'une part, quand $\al<2(1+\mu^{-1})$, les plus grandes valeurs propres ont pour ordre  $N^{\f{1+\mu}{\al}}$, sont asymptotiquement distribu\'ees selon un processus de Poisson et les vecteurs propres associ\'es sont essentiellement port\'es par deux coordonn\'ees (ce ph\'enom\`ene a d\'ej\`a \'et\'e remarqu\'e pour des matrices pleines par  Soshnikov dans \cite{sosh04,sosh06}    quand $\al<2$, et par Auffinger \emph{et al} dans \cite{ABP09} quand $\al<4$). D'autre part, quand $\al>2(1+\mu^{-1})$, les plus grandes valeurs propres ont pour ordre $N^{\f{\mu}{2}}$ et la plupart des vecteurs propres de la matrice sont d\'elocalis\'es, \ie approximativement uniform\'ement distribu\'es sur leurs $N$ coordonn\'ees.
\en{abstract}


\section*{Introduction}
Recently some growing interest has been laid on the understanding of the asymptotic behavior of  both eigenvalues and eigenvectors of random matrices in the large size limit. For Wigner random matrices, that is $N\times N$ Hermitian or real symmetric random matrices with i.i.d. entries (modulo the symmetry assumption), the large-$N$-asymptotic behavior is now well understood, provided the distribution of the entries has sub-exponential decay (or at least a large enough number of moments). It is indeed known from the works of Erd\"os-Schlein-Yau, Tao-Vu  and Knowles-Yin (\cite{ESY1,ESY2,ESY3,TV1,TV2,KY} and see also references therein) that :\\
-- eigenvalues are very close to their theoretical prediction given by well-chosen quantiles of the semi-circle distribution (the proof is based on a \emph{strong semi-circle law}). This also yields universality of local statistics in the bulk and at the edge under some appropriate moment assumptions (see Erd\"os \cite{EReview} e.g. for a review of the recent results).\\
-- eigenvectors are fully delocalized in the following sense. The \emph{localization length}, $L$, of an eigenvector $\bv$ is the typical number of coordinates bearing most of its $\ell^2$ norm. Then it is proved that with very ``high probability" there does not exist an   eigenvector with localization length $L\ll N$. Or roughly speaking all coordinates are in the order of $N^{-1/2}$.\\

In this article, we want to fill in the gap of understanding the role of moments in the delocalization properties of eigenvectors. We will be interested in a model of random matrices that we believe to be quite rich, namely random band matrices with heavy-tailed entries. 

More precisely, the matrices under consideration in this paper are Hermitian random matrices with at most $N^{\mu}$ non zero entries per row. In other words, we force some of the entries of a Wigner matrix to be zero. 
This model is believed to be more complicated than Wigner ensembles due to the fact that there is no reference ensemble: there does not exist a ``simple" band random matrix  ensemble for which eigenvalue/eigenvector statistics can be explicitly computed as for the GUE/GOE in Wigner matrices. Thus usual comparison methods (four moments theorem, Green function comparison method) cannot be used directly in this setting.\\
 Such a model is also believed to exhibit a phase transition, depending on $\mu$. On a physical level of rigor, Fyodorov and  Mirlin \cite{FM} e.g. have explained that for Gaussian entries,  the localization length of a typical eigenvector in the bulk of the spectrum shall be of order $ {L=O(N^{\min(2\mu, 1)})}$ so that eigenvectors should be localized (resp. delocalized or extended) if $\mu <1/2$ (resp. $>1/2$). The only rigorous result in the direction of localization is by Schenker  \cite{Schenker}. Therein it is proved that ${L\ll N^{\min(8\mu,1)}}$ for Gaussian band matrices. On the other hand, delocalization in the bulk is proved by Erd\"os, Knowles,  Yau and  Yin \cite{EKYY} when $\mu>4/5$. In both regimes, it is known from Erd\"os and  Knowles \cite{EK2011CMP, ErdKnow} that typically $L\geq N^{\min(7\mu/6,1)}$ for a certain class of random band matrices (with sub-exponential tails and symmetric distribution).
We refer the reader to Spencer \cite{Spencer} and Erd\"os,  Schlein and  Yau \cite{ESY2} for a more detailed discussion on the localized/delocalized regime. Regarding the edges of the spectrum, much less is known about the  typical localization length of the associated eigenvectors. The authors are not aware of a proof that eigenvectors at the edge are
fully delocalized. However, Sodin's statement \cite{Sodin} combined with
Erd\"os-Knowles-Yau-Yin's results \cite{EKYY} suggest that this should be true when $\mu>5/6.$

We will also allow the i.i.d. non zero entries  to admit only a finite number of moments (which can actually be zero). 
Allowing heavy-tailed entries allows some more localization, especially at the edge of the spectrum, as we can infer from Wigner   matrices. This is discussed in particular in the seminal paper by Cizeau and  Bouchaud \cite{CB}. It is known that the limiting spectral measure of such Wigner matrices is the semi-circle distribution provided that the variance of the entries is finite (otherwise another limiting distribution has been  identified by Guionnet and Ben Arous  \cite{GBAG}).  Regarding eigenvectors, it was shown by  Soshnikov \cite{sosh04,sosh06} and Auffinger,  Ben Arous and  P\'ech\'e \cite{ABP09} that eigenvectors associated to the largest eigenvalues have a localization length of order $1$ if the entries do not admit a finite fourth moment. The localization length is not so clear in the bulk but some progress has been obtained by Bordenave and Guionnet \cite{GB}. However it is commonly believed that the fourth moment shall be a threshold for the localization   of eigenvectors at the edge of the spectrum of full Wigner matrices. For band matrices when the bandwidth is negligible w.r.t. the size $N$ of the matrix, no such threshold has been intuited. This is also a gap we intend to fill in here.\\

Specifically, we prove the following  phase transition, occurring when $\al=2(1+\mu^{-1})$, \ie when $\f{1+\mu}{\al}=\f{\mu}{2}$ (note that $N^{\f{1+\mu}{\al}}$ is always the order of the largest entries of the matrix, while $N^{\mu/2}$ is the order of the bulk of the spectrum when $\al>2$). On the one hand, when  $\al<2(1+\mu^{-1})$, the largest entries of the matrix give rise to isolated eigenvalues with order $N^{\f{1+\mu}{\al}}$ and  eigenvectors essentially carried  by two coordinates. 
 This phenomenon has already been noted  for full matrices ($\mu=1$)  by Soshnikov in \cite{sosh04,sosh06} when $\al<2
$, and by Auffinger \emph{et al} in \cite{ABP09} when $\al<4$. On the other hand, when  $\al>2(1+\mu^{-1})$,  we have  $N^{\f{1+\mu}{\al}}\ll N^{\mu/2}$, so largest entries no longer play any  {specific} role and  the matrix is from this point of view like a matrix with non heavy tailed entries. This is why  the largest eigenvalues have order $N^{\f{\mu}{2}}$ and most eigenvectors of the matrix are delocalized, \ie approximately uniformly distributed on their $N$ coordinates.\\

The paper is organized as follows. In Section \re{section_results_AD}, we state our two main theorems : Theorem \re{Subcritical_case} is the localization result mentioned above about the extreme eigenvalues and eigenvectors in the case   $\al<2(1+\mu^{-1})$ and Theorem \re{Supercritical_case} is the delocalization result mentioned above about the extreme eigenvalues of the matrix  and most of its eigenvectors in the case   $\al>2(1+\mu^{-1})$.  Sections \re{StuckInsideOfMobile_AD1012}, \re{StuckInsideOfMobile_AD10122} and \re{StuckInsideOfMobile_AD10123} are devoted to the proofs of these results and the appendix is devoted to the proof of several technical results, including Theorem \re{147093h48}, a general result whose idea goes back to  papers of Soshnikov  about the surprising phenomenon that certain Hermitian matrices have  approximately equal largest eigenvalues and largest entries.\\

\noindent{\bf Ackowledgment} The authors thank an anonymous referee for very useful comments, which improved the preliminary version of this article.

\noindent{\bf Notation.} For any functions (or sequences) $f,g$, we write  $f(x)\sim g(x)$ (resp. $f(x)\sL g(x)$) when $f(x)/g(x)\lto 1$ (resp. $f(x)/g(x)$ is  \emph{slowly varying})   as $x\to+\infty$. 
We denote by  $\|\bv\|$   the $\ell^2$-norm of $\bv\in\C^N$ and by $\|A\|$   the $\ell^2\to\ell^2$ operator norm of a matrix $A$. When $A$ is normal, then $\|A\|= \rho(A)$ where $\rho(A)$ is the spectral radius of $A$ and we use equivalently both  notations.

An event $A_N$ depending on a parameter $N$ is said to hold  \emph{with exponentially high probability} (abbreviated w.e.h.p. in the sequel) if $\mathbb{P}(A_N)\ge 1-e^{-CN^\tta}$ for some $C,\tta>0$.

\section{The results}\la{section_results_AD}

Let us fix two exponents $\al>0$ and $\mu\in (0,1]$ and let, for each $N$,  $A_N=[a_{ij}]_{i,j=1}^N$   be a real symmetric (or complex Hermitian) random  matrix satisfying the following assumptions:    

{i) For all $i$'s in $\{1, \ld, N\}$ except possibly $o(N)$ of them,  \be\la{largeband12712} \frac{\sharp \{j\ste a_{ij}\trm{ is not almost surely zero}\}}{N^\mu}=a_N ,\ee where $a_N \to 1$ as $N \to \infty.$} {For the other $i$'s, \eqre{largeband12712} still holds, but with $\le$ instead of $=$.} {Thus $A_N$ is a sparse matrix  when $\mu<1$. We denote by $B(N)$ the set of its non a.s. zero entries: 
$$B(N):=\{ (i,j)\ste    a_{ij}\trm{ is not almost surely zero}\}$$ and set $d_N:= a_N N^{\mu}\sim \sharp B(N)/N$.}  

{ii)  The entries $a_{ij}$, $(i,j)\in B(N)$,  are i.i.d. modulo the symmetry assumption  and  \st for a slowly varying function $L$ not depending on $N$,   \be\la{estimate_queue_140712}G(x):=\p(|a_{ij}| \ge x)=L(x) x^{-\al}.\ee}    
If $\al\ge 1+\mu^{-1}$, we also suppose that the $a_{ij}$'s are symmetrically distributed.  This symmetry assumption simplifies the exposition of arguments but can be relaxed (we briefly indicate this possible extension in Remark \ref{rmk: extension} below). At last, if $\al>2(1+\mu^{-1})$, the second moment of the non identically zero entries of $A_N$ is equal to one.  
Note that for all fixed $i,j$, $a_{ij}$ might  depend on $N$ (think for example of the case where $A_N$ is a band matrix), hence should be denoted by $a_{ij}{(N)}$.   

The standard example of matrices satisfying \eqre{largeband12712} is given by \emph{band matrices}, \ie matrices with entries $a_{ij}$ \st $a_{ij}=0$ when  
 $|i-j|> N^\mu/2$. Another very close example is the one of \emph{cyclic band matrices}, \ie matrices with entries $a_{ij}$ \st $a_{ij}=0$ when  
$ |i-j|> N^\mu/2$ and $|i-j|>N-N^\mu/2$.

We denote by $\lam_1\ge \lam_2\ge\cd$ the eigenvalues of $A_N$ (they depend implicitly on $N$) and we choose some unit associated eigenvectors $$\bv_1, \bv_2, \ld$$  Let us also introduce a set of pairs of indices $(i_1\le j_1), (i_2\le j_2),\ld$ \st for all $k$, $|a_{i_kj_k}|$ is the $k$th largest entry, in absolute value, of $A_N$. Let $\tta_k\in \R$ \st $a_{i_kj_k}=|a_{i_kj_k}|e^{2i\tta_k}$.  The eigenvectors $\bv_1, \bv_2, \ld$ are chosen \st for each $k$, $$e^{-i\tta_k}\lan \bv_k, \bfe_{i_k}\ran\ge 0,$$ with $\bfe_1, \ld, \bfe_N$ the vectors of the canonical basis.

As we shall see in the two following theorems, the asymptotic behavior of both the largest eigenvalues of $A_N$ and their associated eigenvectors  exhibit a phase transition with threshold  $\al=2(1+\mu^{-1})$.  \beg{Th}[Subcritical case]\la{Subcritical_case} Let us suppose that $\al<2(1+\mu^{-1})$. Then for each fixed $k\ge 1$, we have the convergences in probability, as $N\to\infty$,  \be\la{371217h29}\f{\lam_k}{|a_{i_kj_k}|}\lto 1\ee and \be\la{371217h210}\qquad \qquad\qquad\qquad\bv_k-\frac{1}{\sqrt 2}\left (e^{i\tta_k} \bfe_{i_k}+e^{-i\tta_k}\bfe_{j_k}\right )  \lto 0 \qquad\trm{ (for the $\ell^2$-norm)}.\ee
As a consequence of \eqre{371217h29}, for $b_N$ the sequence defined by  \eqre{defbn19312} below,  the random point process $$\sum_{k\ste |a_{i_kj_k}|>0}\del_{{\lam_k}/{b_N}}$$ converges in law to the law of a Poisson point process on $(0, +\infty)$ with intensity measure  $\f{\al}{x^{\al+1}}\ud x.$
\en{Th}

The sequence $b_N$ of the theorem is   defined by       \be\la{defbn19312}b_N:= \inf\lf\{ x\ge 0\ste G(x)\le \ff{\#\{\trm{non identically zero  independent entries of $A_N$}\}}\ri\},\ee where $G(x)$ is defined by \eqre{estimate_queue_140712}.   It can easily be deduced from  \eqre{largeband12712} and \eqre{estimate_queue_140712} that \be\la{14712eqb_N}b_N\sL N^{\f{1+\mu}{\al}}.\ee
Roughly speaking, this theorem says that when $\al<2(1+\mu^{-1})$, the largest eigenvalues of $A_N$ have  order $N^{\f{1+\mu}{\al}}$, but no fixed limit when divided by $N^{\f{1+\mu}{\al}}$, because the limiting object is a Poisson process. Moreover, the corresponding eigenvectors are essentially supported by two components. As we shall see in the following theorem, the case $\al>2(1+\mu^{-1})$ is deeply different: in this case, the largest eigenvalues of $A_N$ have order $N^{\f{\mu}{2}}$ and tend to $2$ when divided by $N^{\f{\mu}{2}}$, whereas the eigenvectors are much more delocalized, \ie supported by a large number of components. 

To be more precise,  we use the following Definition 7.1 from  Erd\"os, Schlein and Yau  \cite{ESY2}.
\begin{Def} Let $L$ be a positive integer and $\eta\in (0,1]$ be given. A unit vector $\bv=(v_1, \ld, v_N)\in \C^N$ is said to be   \emph{$(L, \eta)$-localized}  if    there exists a set $S\subset\{1,†\ld, N\}$  \st $|S|=L$ and $\sum_{j\in S^c}|v_j|^2\le \eta$. 
\end{Def}

We shall also use the following slightly modified version of the above definition.  
\begin{Def} Let $L$ be a positive integer and $\eta\in (0,1]$ be given. A unit vector $\bv=(v_1, \ld, v_N)\in \C^N$ is said to be    \emph{$(L, \eta)$-successively localized} if    there exists   a set $S$ which is an interval of the set $\{1,†\ld, N\}$  endowed with the cyclic order \st $|S|=L$ and $\sum_{j\in S^c}|v_j|^2\le \eta$. 
\end{Def}
\beg{rmk}{\rm The larger $L$ and $\eta$, the stronger a statement of the type \emph{``There is non $(L, \eta)$-localized eigenvector"} is.}\en{rmk}

\beg{Th}[Supercritical case]\la{Supercritical_case} Let us suppose that $\al>2(1+\mu^{-1})$ {and that the $a_{ij}$'s have variance one}. Then for each fixed $k\ge 1$,  as $N\to\infty$,  we have the   convergence in probability \be\la{371217h29deloc}\f{\lam_k}{N^{\f{\mu}{2}}}\lto 2.\ee Moreover, for $L:=\lfloor N^c\rfloor$, with  $c$ \st \be\la{condc23612bisLoc}c<\f{2}{5}\mu\f{\al-2}{\al-1} {\qquad \trm{(resp. $c<\mu$)}}\ee 
 for any $\eta_0<1/2$,   we have, as $N\to\infty$,  \be\la{BeatGoesOnAD} \p\lf(\cup_{\eta<\eta_0}\{\exists k, |\lam_k|>\sqrt{2\eta}\rho(A) \trm{ and }\vb_k\trm{ is $(L,\eta)$-localized}\}\ri)\lto 0,\ee
 {(resp. $\ds \p\lf(\cup_{\eta<\eta_0}\{\exists k, |\lam_k|>\sqrt{2\eta}\rho(A) \trm{ and }\vb_k\trm{ is $(L,\eta)$-successively localized}\}\ri)\lto 0$.)}
 \en{Th}
\beg{rmk}{\rm Note that this theorem does not only apply to the edges of the spectrum, as $\eta$ runs from $0$ to $\eta_0$ in  \eqre{BeatGoesOnAD}.}\en{rmk}

 \beg{rmk}{\rm Note that focusing on successively localized vectors, we would need to improve the bound $c<\mu$ in order to get some flavor of the usual threshold of the so-called \emph{Anderson transition}. The localization length $L$ of typical eigenvectors in the bulk is indeed supposed numerically to be in the order of $L\approx N^{2\mu}$ when $\mu<1/2$  for entries with many moments.  At the edge of the spectrum, the authors are not aware of any intuited (even at a physical level of rigor) localization length in the localized regime.}\en{rmk}

To prove both above theorems, we shall also use the following result, which had not appeared at this level of generality yet.

\beg{Th}\la{CVSCL16712}We suppose that the hypotheses \eqre{largeband12712} and \eqre{estimate_queue_140712}  hold with $\al>2$  and that the first and second moments of the non identically zero entries of $A_N$ are respectively equal to $0$ and $1$. Then the empirical spectral measure of $A_N/N^{\f{\mu}{2}}$ converges almost surely to the semi-circle law with support $[-2,2]$.
\en{Th}

\bpr The proof relies on a classical cutt-off and moments method, copying the proof of the convergence to the semi-circle distribution for standard Wigner matrices (see for example \cite[Th. 2.5]{bai-silver-book}).
\epr

\section{A preliminary result: general upper-bound on the moments}\la{StuckInsideOfMobile_AD1012}

 The Hypotheses made on the $a_{ij}$'s are   the ones presented in the beginning of Section \re{section_results_AD}.
 
 \beg{Th}\la{Thmoments22612} Assume that $\al> 2$ and that the $a_{ij}$'s have variance one.   Consider  some positive exponents $\ga,\ga',\ga''$ \st \be\la{1271215h14RH}\f{\mu}{2}\le \ga' \qquad\trm{ and }\qquad \f{\mu}{4}+\ga+\ga''<\ga'\ee and define the truncated matrix  $\hat{A}_N=[a_{ij}\one_{|a_{ij}|\le N^\ga}]_{i,j=1}^N$.
 Then for $s_N\le N^{\ga''}$, there exists a slowly varying function $L_0$ \st 
 $$ \E[\Tr(\hat{A}_N^{2s_N})]\le L_0(N)N^{1+2\gamma} s_N^{-3/2} (2N^{\ga'})^{2s_N}.$$\en{Th}

 The following corollary follows directly from the theorem and from the Chebichev inequality.
\begin{cor}\la{corThmoments22612}
 Under the above hypotheses, there exists a  slowly varying function $L_0$  such that  for any $\ka< 1$ (possibly depending on $N$),   \be\la{14407120444}\p\lf(\|\hat{A}_N\|\ge \ka \ti 2N^{\ga'}\ri)\le \ka^{-2s_N}L_0(N)N^{1+2\gamma} s_N^{-3/2}. \ee  \en{cor}

\beg{rmk}\la{rmkThmoments22612StFlour}{\rm Roughly speaking, this theorem says that  for any $\epsilon>0$,  $$\|\hat{A}_N\|\le(2+\epsilon)N^{\max\{\f{\mu}{2}, \f{\mu}{4}+\ga+\epsilon\}
} \qquad \trm{ for $N\gg1$}.$$}
\en{rmk}
 
\beg{rmk}\la{rmkThmoments22612}{\rm Note that for the theorem and its corollary to be true, one does not really need the size of the matrix to be $N$, but just to be not more than a fixed power of $N$. This remark will allow us to apply the estimate \eqre{14407120444} to submatrices of $A_N$.}
\en{rmk}

 \beg{proofth} Our strategy will be to use the ideas of Soshnikov, well explained  in \cite{soshnisinai98} (see also \cite{soshni99} or \cite{ABP09}).  We shall also  need an estimate on the moments of the truncated variables $\hat{a}_{ij}:=a_{ij}\one_{|a_{ij}|\le N^{\ga}}$. By Lemma \re{lemtruncmomentsHTRV},   there is a slowly varying function $L$ \st 
    for any   $k\ge 0$, for any (non identically null) $a_{ij}$,   \be\la{eqmoments13612} \E[|\ha_{ij}|^{k}]\le \beg{cases} L(N) &\trm{if $k\le \al$,}\\
L(N)N^{\ga(k-\al)}&\trm{if $k>\al$.}\en{cases}\ee

 We have,  suppressing the dependence on $N$ to simplify the notation, $$\Tr \hat{A}^{2s}=\sum_{\substack{1\le i_0, \ld, i_{2s}\le N\\ i_0=i_{2s}}}\ha_{i_0i_1}\cd \ha_{i_{2s-1}i_{2s}}.$$   To any $\bi= (i_0, \ld, i_{2s})$ \st $i_0=i_{2s}$, we associate the non oriented graph $G_\bi:=(V_\bi, E_\bi)$  with vertex set $\{i_0, \ld, i_{2s}\}$ and edges $\{i_{\ell-1}, i_\ell\}$, $1\le \ell\le 2s$ and the closed path $P_\bi=i_0\to i_1\to \cd \to i_{2s}$ on this graph.

Since the $a_{ij}$'s are symmetrically distributed,  each edge of $G_\bi$ has to be visited an even number of times by $P_\bi$ for the contribution of $\bi$ to $\E \Tr \hat{A}^{2s}$ to be nonzero.

 To such a $\bi$, we associate a set $M_\bi$ of $s$ \emph{marked instants} as follows. We read the edges of $ P_\bi$  successively. The instant at which an edge $\{i, j\}$ is read is 
then said to be \emph{marked} if up to that moment (inclusive) the edge $\{i, j\}$ was read an odd number 
of times (note that the instant are counted from $1$ to $2s$, hence the instant where, for example, the edge $i_0\to i_1$ is read is the instant $1$). Other instants are said to be \emph{unmarked}.  Since each edge of $G_\bi$ is visited an even number of times by $P_\bi$, it is clear that $$\# M_\bi=s.$$  Note that at any moment of time $t$, the number of marked instants up to time $t$ is greater (or equal) than the number of unmarked instants. Thus one can associate to the sequence of marked/unmarked instants a unique Dyck path, that is a trajectory in $(\mathbb{Z}^+)^2$, started at $(0,0)$ at time $0$ and ending at $(2s,0)$ at time $2s$  with possible steps $(1,\pm 1)$: 
for any $i=1, \ldots, 2s$, step number $i$ in the trajectory is $(1,1)$ \ssi the instant $i$ is marked. 
 
Now, for each $0\le k\le s$, we define $N_\bi(k)$ to be the set of $i$'s in $\{1,†\ld, N\}$   occurring exactly $k$ times as the current vertex of a marked instant and $n_k$ be its cardinality. Let the family $(n_0, \ld, n_s)$ be called the \emph{type} of $\bi$. 
Note that we have \be\la{mn12312} \sum_{k=0}^{s} n_k=N \qquad\trm{ and }\qquad \sum_{k=0}^{s} k n_k=s.\ee

Let us now count the number of $\bi$'s with fixed type $(n_0, \ld, n_s)$ (where the $n_i$'s satisfy \eqre{mn12312}). To define such a $\bi$, one first has to choose the set $M_\bi$ of marked instants: there are as many possibilities as Dick paths, \ie the Catalan number $$C_s:=\ff{s+1}\binom{2s}{s}.$$
Then one has to choose an unlabelled partition of $M_\bi$ defined by the fact that two marked instants are in the same class if and only if the path $P_\bi$ is at the same vertex of $G_\bi$ at both of these instants. Such a partition is only required to have $n_k$ blocks of cardinality $k$ for each $k=1,†\ld, s$. Hence there are $$ \f{s !}{\prod_{k=1}^s (k!)^{n_k}}\ti\ff{\prod_{k=1}^s n_k!}$$ possibilities (the first factor counting the labelled partitions and the second one ``delabelling'').
At this point, one has to choose the vertices of $G_\bi$. For $i_0$, there are $N$ possibilities. For each other vertex, there are at most $d_N$ possibilities. There are at most $n_1+\cd+n_s$ other vertices. Indeed, except possibly $i_0$, each vertex is occurring a certain number of  times as the current vertex of a marked instant (for example at the first time the vertex is visited by the path $P_\bi$). Hence there are at most $Nd_N^{\,n_1+\cd+n_s}$ possibilities for the choices of the vertices of $G_\bi$. 
There now remains to give an upper-bound on the number of ways to determine vertices at unmarked instants (such vertices will not be new, but still have to be chosen among the before chosen vertices). Soshnikov proved in \cite{soshni99}, Formula (4.3), or \cite{Sos}, first paragraph on p. 6, that this number is not larger than $\prod_{k=2}^s(2k)^{kn_k}$ (the idea is that if $v$ is a vertex arising at $k$ marked instants, the number of ways to determine the endpoint of an edge starting from $v$ at an unmarked instant is at most $2k$).

To sum up, the number of $\bi$'s with fixed type $(n_0, \ld, n_s)$ is at most $$C_s\f{s !}{\prod_{k=1}^s (k!)^{n_k}}\ti\ff{\prod_{k=1}^s n_k!}\ti
Nd_N^{\, n_1+\cd+n_s}\ti\prod_{k=2}^s(2k)^{kn_k}.$$

 Let us now give an upper bound on the expectation $\E[\ha_{i_0i_1}\cd \ha_{i_{2s-1}i_{2s}}]$ depending on the type of $\bi$. For $i,j\in V_\bi$, let $\ovl{ij}$ denote the edge $\{i,j\}$ of  $G_\bi$ (this edge is unoriented, so $\ovl{ij}=\ovl{ji}$) and let $k({\ovl{ij}})$ denote the half of number of times that this edge is visited by $P_\bi$, i.e. the number of marked instants along edge $\ovl{ij}$. We also introduce $k_{i;\ovl{ij}}$ to be the number of times that the vertex $i$ is marked along the edge $\ovl{ij}$. Clearly, $$k({\ovl{ij}})=k_{i;\ovl{ij}}+k_{j;\ovl{ij}}\quad\trm{ and }\quad \op{type}(i)=\sum_jk_{i;\ovl{ij}}.$$
 We know, by \eqre{eqmoments13612}, that for a certain slowly varying sequence $L(N)$ (that can change at every line)
\begin{eqnarray*}&  \E[\ha_{i_0i_1}\cd \ha_{i_{2s-1}i_{2s}}]& = \prod_{\substack{e\in E_\bi,\\ k(e)\ge 2}}L(N)\prod_{\substack{e\in E_\bi,\\ k(e)\ge\al }}N^{\ga(2k_e-\al)} \\
&&\le  \prod_{\substack{e\in E_\bi,\\ k(e)\ge 2}} L(N)\prod_{\substack{e\in E_\bi,\\ k(e)\ge 2 }}N^{\ga(2k_e-2)} =  L(N)^{2E} N^{-2\ga E}\prod_{\substack{e\in E_\bi,\\ k(e)\ge 2 }}N^{2\ga k_e},
\end{eqnarray*}
 where $E$ denotes the number of edges $e$ \st $k(e)\ge 2$.
Let us now enumerate the edges via their extremities.   
 Then \beq\sum_{\substack{e\in E_\bi,\\ k(e)\ge 2 }} 
 k(e)&=&\#\{\trm{marked instants along edges $e$ \st $k(e)\ge 2$}\}\\
 &=& \sum_{\substack{(v,w)\in V_\bi^2  \\ k({\ovl{vw}})\ge 2}}k_{v; \ovl{vw}}
 = \sum_{\substack{(v,w)\in V_\bi^2 \\ k({\ovl{vw}})\ge 2,\op{type}(v)=1}}k_{v; \ovl{vw}}\;\;+ \sum_{\substack{(v,w)\in V_\bi^2 \\ k({\ovl{vw}})\ge 2,\op{type}(v)\ge 2}}k_{v; \ovl{vw}}.\eeq 
 Let us now use the fact, well known from \cite{soshni99,soshnipeche07,ABP09} that if an edge $\ovl{vw}$ is visited at least $4$ times by the path $P_\bi$, then at least one of $v$ and $w$ have type $\ge 2$, except for the first visited vertex $i_0$. It follows that the first sum above is $\le E+1$ and also that $E\le 1+\sum_{k=2}^s kn_k$. Hence 
 \beq
\sum_{\substack{e\in E_\bi,\\ k(e)\ge 2 }} 
 k(e) &\le & E+1\quad+\sum_{v\in V_\bi, \op{type}(v)\ge 2}\sum_w k_{v;\ovl{vw}}\\
   &=&E+1\quad +\sum_{k=2}^s kn_k \eeq  
 Hence  $  {\E[\ha_{i_0i_1}\cd \ha_{i_{2s-1}i_{2s}}] \le L(N)^{2+2\sum_{k=2}^skn_k}
 N^{2\ga (1+\sum_{k=2}^s kn_k) }},$
 using   $E\leq 1+\sum_{k\geq 2}kn_k.$ 

As a consequence, \begin{eqnarray*}&\E[\Tr(\hat{A}^{2s})]/N^{2s\ga'}& \le  L(N)^2 N^{1+2\ga-2s\ga'}C_ss!\!\!\!\!\!\!\sum_{\substack{n_1, \ld, n_s\\ \trm{s. t. \eqre{mn12312} holds}}} \!\!\!\!{\ff{\prod_{k=1}^s (k!)^{n_k}}\ti\ff{\prod_{k=1}^s n_k!}\ti d_N^{n_1+\cd+n_s}\ti\prod_{k=2}^s(2k)^{kn_k}} \cr
&& \ti N^{2\ga \sum_{k=2}^s kn_k } L(N)^{2\sum_{k=2}^s kn_k}. \end{eqnarray*}
Let  us now use the fact  $s!\le n_1!s^{(s-n_1)}=n_1!s^{\sum_{k=2}^skn_k}\le n_1!N^{\ga'' \sum_{k=2}^skn_k}$,  $d_N\le N^\mu$    and  $N^{-2s\ga'}=N^{-2\ga'\sum_{k=1}^s kn_k}$. 
We get  $$\E[\Tr(\hat{A}^{2s})]/N^{2s\ga'}  \le$$ $$ L(N)^{2} N^{1+2\ga}C_s \sum_{\substack{n_1, \ld, n_s\\ \trm{s. t. \eqre{mn12312} holds}}} {N^{(\mu-2\ga')n_1}}  \prod_{k=2}^s \ff{n_k!} \lf(\f{ L(N)^{2k}N^\mu (2kN^{2(\ga-\ga')+\ga''})^k}{k!}\ri)^{n_k}.
$$ 
But by the hypothesis \eqre{1271215h14RH}, ${\mu-2\ga'}\le 0$, hence the first factor is $\le 1$, so, using the fact that by \eqre{mn12312}, $n_1$ is determined by the other $n_j$'s, we get 
 \beq \E[\Tr(\hat{A}^{2s})]/N^{2s\ga'} &\le & L(N)^{2} N^{1+2\ga} C_s\sum_{ n_2, \ld, n_s\ge 0} \prod_{k=2}^s \ff{n_k!} \lf(\f{ L(N)^{2k}N^\mu (2kN^{2(\ga-\ga')+\ga''})^k}{k!}\ri)^{n_k} 
\\
&\le & L(N)^{2} N^{1+2\ga} C_s\sum_{ n_2, \ld, n_s\ge 0} \prod_{k=2}^s \ff{n_k!} \lf(\f{ L(N)^{2k}N^\mu (2N^{2(\ga-\ga'+\ga'')})^k}{k!}\ri)^{n_k} 
\\
&\le &  L(N)^{2} N^{1+2\ga} C_s \exp\lf(\sum_{k=2}^s \f{ L(N)^{2k}N^\mu (2N^{2(\ga-\ga'+\ga'')})^k}{k!}\ri)\\
&\le &  L(N)^{2} N^{1+2\ga} C_s \exp\lf(\sum_{k=2}^s \f{L(N)^{2k}2^kN^{-\eps k}}{k!}\ri) \eeq   
 with $\eps=-2\lf(\f{\mu}{4}+\ga+\ga''-\ga'\ri)$. By  \eqre{1271215h14RH}, we have $\eps>0$, so that the exponential term stays bounded as $N\to\infty$.  Using the fact that $L(x)^2$ varies slowly and   $C_s\sim 4^s (\pi s^3)^{-1/2}$, we get Theorem \ref{Thmoments22612}.\en{proofth}

\beg{rmk}\label{rmk: extension}{\rm 
In the case where the entries $a_{ij}, 1\leq i,j \leq N$ are not symmetrically distributed, one can prove a similar statement as in 
Theorem \ref{Thmoments22612}. The proof is based on arguments already given in Section 4 of \cite{ABP09} and \cite{soshnipeche07}. One can indeed assume that the truncated entries are centered. Then, the main modification in evaluating $\E[\Tr(\hat{A}^{2s})]/N^{2s\ga'}$ is that one has to take into account the contribution of paths with edges seen an odd number of times. However any such edge is seen at least 3 times, because the entries are centered. It can then be shown that the contribution of such paths is negligible (provided $s$ is small enough as in Theorem \ref{Thmoments22612}), as to each such edge corresponds a vertex of type $>1$.}\en{rmk}

\section{Proof of Theorem \re{Subcritical_case}}\la{StuckInsideOfMobile_AD10122}
We are going to prove Theorem \re{Subcritical_case} as an application of Theorem \re{147093h48} of the appendix, for $c_n:=b_N$, the sequence defined by \eqre{defbn19312}. More precisely, we shall use its ``random versions": the case   $\al<1+\mu^{-1}$ will be a consequence of Corollary \re{tr2131212h13StFlour}, whereas in the case   $1+\mu^{-1}\le\al< 2(1+\mu^{-1})$, we need to truncate the entries, and the conclusion   will follow from Corollary \re{tr2131212h13}.

Hypothesis \eqre{estimate_queue_140712} implies  that the distribution of the non-zero entries is in the max-domain of attraction of the Fr\'echet distribution with exponent $\al$ (see \cite{resnick}, p. 54). By e.g. \cite[Th. 2.3.1]{LLRextremes}, it implies that as $N\to\infty$, the point process \bes\la{1671214h43}\sum_{k\ste |a_{i_kj_k}|>0}\del_{{|a_{i_kj_k}|}/{b_N}}\ees converges in distribution to  a Poisson point process on $(0, +\infty)$ with intensity measure  $\f{\al}{x^{\al+1}}\ud x.$ It explains why the second part of Theorem \re{Subcritical_case} is a consequence of its first part and why Hypothesis 
\eqre{1471210h19} of the corollaries \re{tr2131212h13StFlour} and \re{tr2131212h13}  is satisfied by the $|a_{i_kj_k}|$'s.

Note that by \eqre{estimate_queue_140712},  for any $\tta >0$, for any non indenticaly null $a_{ij}$, \be\la{eq_tail_b}\p(|a_{ij}|>b_N^\tta)\sL b_N^{-\al \tta}\sL N^{-\tta(1+\mu)}.\ee
The following claim (valid without any assumption on $\al$)  is a direct consequence of \eqre{eq_tail_b} and of the union bound. 

\beg{claim}\la{1670912h192}For any $\eta>0$,  with \pro going to one as $N\to\infty$, we have: \bgt\ite[a)] no row of $A_N$ has two entries larger, in absolute value, than $b_N^{\f{1+2\mu}{2(1+\mu)}+\eta}$ (the exponent $\f{1+2\mu}{2(1+\mu)}$ increases from $\ff{2}$ to $\f{3}{4}$ as $\mu$ increases from $0$ to $1$),\\
\ite[b)] the matrix $A_N$ has no diagonal entry larger, in absolute value, than $b_N^{\ff{1+\mu}+\eta}$.
\ent\en{claim}

So  for any positive $\del,\del'$,  parts  (b.i) and (b.ii) of the random version of Hypothesis  \re{h531215h282} are satisfied with \be\la{defbn193122}\ka:= \f{1+2\mu}{2(1+\mu)}+\del, \qquad\tau:=\ff{1+\mu}+\del'.\ee  Let us now verify  Part (b.iii). 

\subsection{Case where $\al<1+\mu^{-1}$} 
Set  \be\la{1471212h}S:=\sum_{j\ste |a_{ij}|<b_N^{\kappa}}|a_{ij}|\ee (we suppress the dependence in $N$ to simplify notation). We shall prove that there is $\nu<1$ \st   $\p(S_N>b_N^\nu)$ is exponentially small with some bounds that are uniform  on $i$.   The sum $S$ can be rewritten $S=S_1+S_2+S_3 $ as follows : $$S=\sum_{|a_{ij}|\le N^{\f{\mu}{\al}-\eta}}|a_{ij}| +\sum_{ N^{\f{\mu}{\al}-\eta}<|a_{ij}|\le N^{\f{\mu}{\al}+\eta}}|a_{1i}|+
\sum_{ N^{\f{\mu}{\al}+\eta}<|a_{ij}|\le b_N^{\ka}}|a_{ij}| $$
The sums $S_1, S_2, S_3$ can be treated with respectively parts a), c) and d) of Proposition \re{proposumHTV} of the appendix.  The treatment of $S_1$ uses the facts $b_N\sL N^{\f{1+\mu}{\al}}$ and  that $\mu+\f{\mu}{\al}(1-\al)^+<\f{1+\mu}{\al},$ which is always true when $\al\le 1$ and which is a consequence of  $\al<1+\mu^{-1}$ when $\al>1$.

\subsection{Case where $1+\mu^{-1}\le \al<2(1+\mu^{-1})$}   We have seen at \eqre{14712eqb_N} that $b_N\sL N^{\f{1+\mu}{\al}}$. So
to apply Corollary \re{tr2131212h13} for $c_n=b_N$, we have to find a cut-off exponent $\ga$ satisfying both following constraints: \bgt\ite[1)]
for $\ka$ defined by \eqre{defbn193122}, there is $\eps>0$ \st with exponentially high  probability,   \be\la{591023h17}S:=\sum_{j\ste N^\ga<|a_{ij}|\le b_N^\ka}|a_{ij}|\le N^{\f{1+\mu}{\al}-\eps},\ee
\ite[2)] there is $\eps'>0$ \st with probability tending to one, we have:  $$\qquad\qquad\qquad\qquad\|\hat{A}_N\|\le N^{\f{1+\mu}{\al}-\eps'}\qquad\qquad\trm{ (with $\hat{A}_N:=[a_{ij}\one_{|a_{ij}|\le N^\ga}]_{1\le i,j\le N}$)}.$$
\ent 

By Corollary \re{corThmoments22612}, the second condition is satisfied when $\max\{\f{\mu}{2}, \f{\mu}{4}+\ga\}<\f{1+\mu}{\al}$.
As $\al<2(1+\mu^{-1})$, we have $\f{\mu}{2}<\f{1+\mu}{\al}$, so for condition 2) to be verified, one only needs $$\ga<\f{1+\mu}{\al}-\f{\mu}{4}.$$
 To treat the sum $S$ of \eqre{591023h17}, one proceeds as we did to treat the sum $S$ defined at \eqre{1471212h} in the case $\al<1+\mu^{-1}$, except that now,  $$S_1=\sum_{N^\ga<|a_{1j}|\le N^{\f{\mu}{\al}-\eta}}|a_{1j}|$$ and $S_1$ is treated thanks to Part b)  of Proposition \re{proposumHTV}. Indeed, Part b)  of Proposition \re{proposumHTV} implies that  w.e.h.p., $S_1\le N^{\mu-\ga(\al-1)+\eta}$ with $\eta>0$ as small as we need. 
 Hence to fulfill Condition 1),  one needs $\ga$ to satisfy 
 \bes\mu-\ga(\al-1)<\f{1+\mu}{\al},\ees i.e.   $\ga>\f{\mu}{\al}-\ff{\al(\al-1)}.$ To sum up, we need to find a cut-off exponent $\ga$ \st $$\f{\mu}{\al}-\ff{\al(\al-1)}<\ga< \f{1+\mu}{\al}-\f{\mu}{4}.$$Hence to conclude, it suffices to remark that $   \al  <  2(1+\mu^{-1})$ implies  \be\la{encalconsec} \f{\mu}{\al}-\ff{\al(\al-1)}<\f{1+\mu}{\al}-\f{\mu}{4}.\ee

\section{Proof of Theorem  \re{Supercritical_case}}\la{StuckInsideOfMobile_AD10123}

\subsection{Eigenvalues}Let us first prove the part about the eigenvalues, \ie Equation \eqre{371217h29deloc}. 
First, by Theorem \re{CVSCL16712}, for any fixed $k\ge 1$, we have $$\liminf\f{\lam_k}{N^{\f{\mu}{2}}}\ge 2.$$
Let us now prove that \be\la{1670911h25}
\limsup\f{\lam_k}{N^{\f{\mu}{2}}}\le 2.
\ee
To do that, we will prove that one can find a cut-off exponent $\ga$ \st for $\hat{A}_N:=[a_{ij}\one_{|a_{ij}|\le N^\ga}]_{1\le i,j\le N}$, we have \be\la{1670911h38}
\limsup\f{\|\hat{A}_N\|}{N^{\f{\mu}{2}}}\le 2 \qquad\trm{ and }\qquad  {\|A_N-\hat{A}_N\|}=o(N^{\f{\mu}{2}}).
\ee
To treat the first part of \eqre{1670911h38}, we apply Corollary \re{corThmoments22612} with $\ga'=\f{\mu}{2}$ and $\ga''>0$ such that $$
\f{\mu}{4}+\ga+\ga''<\ga'.$$For such a $\ga''$ to exist, the constraint on $\ga$ is that $\ga<\f{\mu}{4}$.	
To treat the second part of \eqre{1670911h38}, we use the following claim and the fact that \be\la{1671221j}\|A_N-\hat{A}_N\|=\sup_{\lam\trm{ eig. of }A_N-\hat{A}_N}|\lam|\le \|A_N-\hat{A}_N\|_{\ell^\infty\to\ell^\infty}\le \max_i\sum_{j\ste |a_{ij}|>N^\ga}|a_{ij}|.\ee
 \beg{claim}\la{1671216h20}
Under the hypothesis that  $\al>2(1+\mu^{-1})$, for any $\ga>\f{\mu}{2(\al-1)}$,  there is $\eta>0$ \st with \pro tending to one, we have:  $$\max_i\sum_{j\ste |a_{ij}|>N^\ga}|a_{ij}|\le N^{\f{\mu}{2}-\eta}.$$
\en{claim}
Let us conclude the proof of the eigenvalues part of Theorem  \re{Supercritical_case} before proving the claim. All we need is to find a cut-off exponent $\ga$ \st $$\f{\mu}{2(\al-1)}<\ga<\f{\mu}{4}.$$ The existence of such a $\ga$ is equivalent to the fact that $\al-1>2$, which is true because $\al>2(1+\mu^{-1})\ge 4$.

{\noindent {\it Proof of the claim. }} Let $S(i)$ be the sum in the statement.  By \eqre{estimate_queue_140712}, it is easy to see that for any $\tta>\f{1+\mu}{\al}$, with \pro tending to one, we have $$\max_{ij}|a_{ij}|\le N^{\tta}.$$
Hence by Part a) of Claim \re{1670912h192} (using the fact that $b_N\sL N^{\f{1+\mu}{\al}}$), for such a $\tta$, with \pro tending to one, for all $i$, $$S(i)\le  \sum_{j\ste N^\ga<|a_{ij}|\le N^\tta} |a_{ij}|= \sum_{j\ste N^\ga<|a_{ij}|\le N^{\f{\mu}{\al}}} |a_{ij}|+\sum_{j\ste N^{\f{\mu}{\al}}<|a_{ij}|\le N^{\f{1+\mu}{\al}}} |a_{ij}|+\sum_{j\ste N^{\f{1+\mu}{\al}}<|a_{ij}|\le N^\tta} |a_{ij}|.$$
 Using respectively parts b), c), d) of Proposition \re{proposumHTV}, for any $\epsilon>0$, we have, w.e.h.p.,
$$\sum_{j\ste N^\ga<|a_{ij}|\le N^{\f{\mu}{\al}}} |a_{ij}|\le N^{\mu-\ga(\al-1)+\epsilon},\qquad \sum_{j\ste N^{\f{\mu}{\al}}<|a_{ij}|\le N^{\f{1+\mu}{\al}}} |a_{ij}|\le N^{\f{1+\mu}{\al}+\epsilon}$$ and $$\sum_{j\ste N^{\f{1+\mu}{\al}}<|a_{ij}|\le N^\tta} |a_{ij}|\le N^{\tta+\epsilon}. $$
Hence as $\tta>\f{1+\mu}{\al}$,  for any $\phi>\max\{\mu-\ga(\al-1),\tta\}$, we have, w.e.h.p., uniformly on $i$, $$\sum_{j\ste N^\ga<|a_{ij}|\le N^\tta} |a_{ij}|\le N^\phi.$$Now, to conclude the proof of the claim, it suffices to notice that the hypotheses  $\ga>\f{\mu}{2(\al-1)}$ and  $\al>2(1+\mu^{-1})$  are respectively equivalent to 
   $\mu-\ga(\al-1)<\f{\mu}{2}$ and   $\f{1+\mu}{\al}<\f{\mu}{2}$, so that one can  find some exponents  $\tta,\phi$ satisfying 
$$\f{1+\mu}{\al}<\tta \qquad \trm{ and }\qquad \max\{\mu-\ga(\al-1),\tta\}<\phi<\f{\mu}{2}.$${\hfill $\square$\\}
 
\subsection{Eigenvectors}

We shall first prove the following lemma. Let us recall that a \emph{principal submatrix} of a matrix $H=[x_{ij}]_{1\le i,j\le N}$ is a matrix of the type 
   $H=[x_{j_kj_\ell}]_{1\le k,\ell\le L}$, where $1\le L\le N$ and $1\le j_1<\cdots<j_L\le N$.  The submatrix will be said to be \emph{successively extracted} if  the indices $j_1, \ld, j_L$ form an interval of the set $\{1, \ld, N\}$ endowed with the cyclic order.

  \beg{lem}\la{22101219h15}  Let $H$ be a Hermitian matrix  and $\rho_L(H)$ (resp. $\rho^{succ}_L(H)$)  be the maximum spectral radius of its $L\ti L$ principal (resp. principal successively extracted) submatrices. Let $\lam$ be an eigenvalue of $H$ and $\bv$ an associated unit eigenvector.     

If $\bv$ is $(L, \eta)$-localized, then 
$|\lam|\le \f{\rho_L(H)+\sqrt{\eta}\rho(H)}{\sqrt{1-\eta}}.$ 

If $\bv$ is $(L, \eta)$-successively localized, then $|\lam|\le \f{\rho^{succ}_L(H)+\sqrt{\eta}\rho(H)}{\sqrt{1-\eta}}$.
\en{lem}

\bpr Let $j_1< \cd< j_L$ be  indices \st  $\sum_{\ell=1}^L|v_{j_\ell}|^2\ge 1-\eta$ and let $P$ be the orthogonal projection onto the subspace generated by the vectors $\bfe_{j_1}, \ld, \bfe_{j_L}$ (the $\bfe_j$'s are the vectors of the canonical basis). We have 
$$\lam P\bv=PH\bv=PHP\bv+PH(1-P)\bv.$$
Then the conclusion follows directly from the following $$|\lam|\ti \sqrt{1-\eta}\le |\lam|\ti \|P\bv\|\le \rho(PHP)+\rho(H)\|(1-P)\bv\|\le \rho(PHP)+\sqrt{\eta}\rho(H).$$\epr

        \beg{claim}Let us suppose   that $\al>2(1+\mu^{-1})$.
\bgt\ite[a)]
     Let us fix $c$ \st \be\la{condc23612bis}c<\f{2}{5}\mu\f{\al-2}{\al-1}.\ee Then there is $\eps>0$ \st w.e.h.p.,  the following holds:\\\begin{quote}For any $\lfloor N^c\rfloor\ti\lfloor N^c\rfloor$ principal submatrix $B$ of $A_N$,   $\|B\|\le N^{\f{\mu}{2}-\eps}$.\\ 
    \end{quote}
\ite[b)]Let us fix $c$ \st \be\la{condc23612bisterter}c<\mu.\ee Then there is $\eps>0$ \st w.e.h.p.,  the following holds:\\ \begin{quote}For any $\lfloor N^c\rfloor\ti\lfloor N^c\rfloor$ successively extracted principal submatrix $B$ of $A_N$, we have  $\|B\|\le N^{\f{\mu}{2}-\eps}$.\\ 
    \end{quote}
\ent
    \en{claim}
  
Before proving the claim, let us conclude to the proof of the eigenvectors part of Theorem  \re{Supercritical_case}. We know that $\rho(A)\sim 2N^{\f{\mu}{2}}$ and that there is $\eps>0$ \st  with  probability tending to one,  $\rho_L(A)$ (resp. $\rho_L^{\trm{succ}}(A)$) is bounded from above by $N^{\f{\mu}{2}-\eps}$.  Since $\ff{\sqrt{1-\eta}}<\ff{\sqrt{1-\eta_0}}<\sqrt{2}$,   the two assumptions $|\lam_k|>\sqrt{2\eta}\rho(A)$ and $\vb_k\trm{ is $(L,\eta)$-localized}$ are then incompatible by 
 Lemma \re{22101219h15}. The case of successively extracted eigenvectors is handled similarly.

 {\noindent {\it Proof of the claim. }}  We shall treat a) and b) in the same time.   Let us first note that Equation \eqre{condc23612bis} (resp. Equation \eqre{condc23612bisterter}) is equivalent to (resp. implies that) $$\f{c}{4}+\f{\mu}{2(\al-1)}+c<\f{\mu}{2}\qquad\trm{(resp. $\f{c}{4}+\f{\mu}{2(\al-1)}<\f{\mu}{2}$)}.$$ Hence one can choose some positive exponents $\eps,\ga,\ga',\ga''$ \st $\gamma' \geq c/2$ and 
   \be\la{condgamma23612}\ga>\f{\mu}{2(\al-1)}, \quad \ga''>c \quad\trm{ (resp. $\ga''>0$) }\quad\trm{ and }\quad \f{c}{4}+\ga+\ga''<\ga'<\f{\mu}{2}-\eps.\ee 
   
   Any submatrix $B=[a_{j_k j_\ell}]_{1\le k,\ell\le\lfloor N^c\rfloor} $ can be written $B=\hat{B}+(B-\hat{B})$, with $\hat{B}:=[a_{j_k j_\ell}\one_{|a_{j_k j_\ell}|\le N^\ga}]_{k,\ell}$. We know (see e.g. \eqre{1671221j})  that  independently of the choice of the $j_k$'s $$\|B-\hat{B}\|\le \max_{1\le i\le N} \sum_{j\trm{ s.t. }|a_{ij}|>N^\ga}|a_{ij}|.$$
Hence by Claim \re{1671216h20},  the condition $\ga>\f{\mu}{2(\al-1)}$ of Equation \eqre{condgamma23612} ensures us that for a certain  $\eta>0$,  with   probability tending to one,   independently of the choice of the $j_k$'s $\|B-\hat{B}\|\le N^{\f{\mu}{2}-\eta}$. Hence  one can focus on $\hat{B}$.
  
  Let us now apply  Corollary \re{corThmoments22612} and Remark \re{rmkThmoments22612}. We get that for any choice of $j_1, \ld, j_{\lfloor N^c\rfloor}$, up to a polynomial factor in the RHT, $$\p(\|\hat{B}\|\ge N^{\f{\mu}{2}-\eps})\le N^{-(\f{\mu}{2}-\eps-\ga')2\lfloor N^{\ga''}\rfloor}.$$ But there are at most $N^{N^c}$ (resp. $N$) ways to choose the indices  $j_1, \ld, j_{\lfloor N^c\rfloor}$ of the rows of the submatrix $B$ (resp. of the successively extracted submatrix $B$). Hence the probability that  $\|\hat{B}\|\ge N^{\f{\mu}{2}-\eps}$ for at least one of these choices  is $\le N^{-(\f{\mu}{2}-\eps-\ga')2\lfloor N^{\ga''}\rfloor+N^c}$ (resp. $\le N^{-(\f{\mu}{2}-\eps-\ga')2\lfloor N^{\ga''}\rfloor+1}$). Since by \eqre{condgamma23612}, $\ga''>c$ (resp. $\ga''>0$) and  $\f{\mu}{2}-\eps-\ga'>0$, the conclusion follows.
  {\hfill $\square$\\}

\section{Appendix}

\subsection{Eigenvalues and eigenvectors under perturbation}
In this section, we state a result about eigenvectors and eigenvalues of perturbed Hermitian matrices. The part about eigenvalues can be found in the literature (see  the books by Bhatia \cite{bhatiamatrixanalysis,bhatiapertmatrices}), but we did not find the part about the eigenvectors in the literature. 

\beg{propo}\la{pr_crucial}
Let $H$ be a Hermitian matrix and $\bv$ be a unit vector \st for a certain $\lam\in \R$, $$H\bv=\lam \bv+\eps \mathbf{w},$$ with $\mathbf{w}$ a unit vector \st $\mathbf{w}\perp \bv$ and $\eps>0$.\bgt\ite[a)] Then $H$ has an eigenvalue $\lam_\eps$ in the ball $\bar{B}(\lam, \eps)$.\\
\ite[b)] Suppose moreover that $H$ has only one eigenvalue (counted with multiplicity) in $\bar{B}(\lam, \eps)$ and that all other eigenvalues are at distance at least $d>\eps$ of $\lam$. Then for $\bv_\eps$ a unit eigenvector associated to $\lam_\eps$, we have $$\|\bv_\eps-P_\bv(\bv_\eps)\|\le \f{2\eps}{d-\eps},$$ where $P_\bv$ denotes the orthogonal projection onto $\op{Span}(\bv)$.\ent  \en{propo}

\bpr Part a) is a simple consequence of perturbation theory (see e.g. Lemma A.2 in \cite{KY2012}).   To prove the second part, note that 
$$\|\bv_\eps-P_\bv(\bv_\eps)\|=\sqrt{1-|\langle \bv_{\eps},\bv \rangle|^2}$$ is symmetric in $\bv,\bv_\eps$. 
Thus we decompose $\bv=\langle \bv_{\eps},\bv \rangle \bv_\eps+\mathbf{r}$ with $\mathbf{r}\perp \bv_\eps$, and prove that \be\la{2406151}\|\mathbf{r}\|\;\le\;  \f{2\eps}{d-\eps}.\ee
We have $$\lam \bv+\eps \mathbf{w}=H\bv=\lam_\eps\langle \bv_{\eps},\bv \rangle \bv_\eps+H\mathbf{r},$$
so that $$\lam\langle \bv_{\eps},\bv \rangle \bv_\eps+\lam\mathbf{r}+\eps \mathbf{w}= \lam_\eps\langle \bv_{\eps},\bv \rangle \bv_\eps+H\mathbf{r},$$
and that $$(\lam-H)\mathbf{r}=(\lam_\eps-\lam)\langle \bv_{\eps},\bv \rangle \bv_\eps-\eps \mathbf{w}.$$
This yields the resut, by considering the norm of $(\lambda-H)$ restricted to the subspace $\bv_\eps^{\perp}$.
\epr

\subsection{Largest eigenvalues vs largest entries of matrices}
In this section, we present a synthetic version of some ideas first appeared in Soshnikov's paper \cite{sosh04}. We also extend these ideas to the eigenvectors level. Theorem \re{147093h48} below gives a sufficient condition for a large deterministic Hermitian matrix to have its $k$th largest eigenvalue approximately equal to its $k$th largest entry in absolute value for all fixed $k$. Note that this is not what happens usually: in some way the large entries need to overwhelm the other entries.\\
We also give sufficient condition so that the corresponding eigenvector is approximately equal to the eigenvector of the symmetric matrix formed by forcing all but this $k$th largest entry to be $0$.  The sufficient condition is, roughly speaking, that the largest entries and their spacings have an order $c_n\gg1$, are sufficiently well spread out in the matrix and that, up to the removing of these largest entries, the sum of the terms of each row of the matrix have order $\ll c_n$. In Corollary \re{tr2131212h13StFlour}, we give the random matrix version of this theorem and in Corollary \re{tr2131212h13}, we explain how one is allowed to first remove a part of the matrix which does not affect the largest entries.

For each $n$, let $H_n$  be an $n\ti n$ deterministic Hermitian matrix with entries $h_{ij}$. 
Let us denote by $\lam_1\ge \lam_2\ge\cd$ the eigenvalues of $H_n$ (they depend implicitly on $n$) and let us choose some unit associated eigenvectors $$\bv_1, \bv_2, \ld$$ Let us also introduce a set of pairs of indices $(i_1\le j_1), (i_2\le j_2),\ld$ \st for all $k$, $|h_{i_kj_k}|$ is the $k$th largest entry, in absolute value, of $H_n$. Let $\tta_k\in \R$ \st $h_{i_kj_k}=|h_{i_kj_k}|e^{2i\tta_k}$ The eigenvectors $\bv_1, \bv_2, \ld$ are chosen \st for each $k$, $e^{-i\tta_k}\lan \bv_k, \bfe_{i_k}\ran\ge 0.$
   We make the following hypotheses.

\beg{hyp}\la{h531215h282}
\bgt
\ite[(a)] There is a sequence    $c_n\lto+\infty$ \st  for any fixed $k$,\\ \bgt \ite[(a.i)]  $c_{n+k}\sim c_n$,\\  \ite[(a.ii)]    $0<\liminf \f{|h_{i_kj_k}|}{c_n}\le \limsup \f{|h_{i_kj_k}|}{c_n}<\infty$ and $\liminf \f{|h_{i_kj_k}|-|h_{i_{k+1}j_{k+1}}|}{c_n}>0$\\\ent
\ite[(b)] There exists  three exponents $ \kappa,\tau, \nu\in (0,1)$  \st for $n$ large enough,\\  \bgt \ite[(b.i)]  no row of $H_n$ has two entries larger, in absolute value, than $c_n^{\kappa}$,\\
\ite[(b.ii)] no diagonal entry of $H_n$ is larger, in absolute value, than $c_n^{\tau}$,\\
\ite[(b.iii)]  for each  $i\in†\{1, \ld, n\}$, $$\sum_{j\ste |h_{ij}|<c_n^{\kappa}}|h_{ij}|\le c_n^{\nu}.$$\ent
\ent 
\en{hyp}

\beg{Th}\la{147093h48}Under Hypothesis \re{h531215h282}, as $n\to\infty$, for any $k\ge 1$ fixed, $$\f{\lam_k}{|h_{i_kj_k}|}\lto 1\qquad\trm{and}\qquad \bv_k-\f{e^{i\tta_k}\bfe_{i_k}+e^{-i\tta_k} \bfe_{j_k}}{\sqrt{2}} \lto 0 \qquad\trm{ (for the $\ell^2$-norm)}.$$
\en{Th} 
Before proving the theorem, let us state its two ``random versions". 
\beg{cor}\la{tr2131212h13StFlour} Suppose now that the matrix $H_n$ is random, that the sequence $b_n$ is  deterministic and satisfies Hypothesis (a.i), replace Hypothesis (a.ii)  by \be\la{1471210h19}\lim_{\eps\to 0} \limsup_{n\to\infty}\p(\f{|h_{i_kj_k}|}{c_n}<\eps)+ \p(\f{|h_{i_kj_k}|}{c_n}>\ff{\eps})+\p(\f{|h_{i_kj_k}|-|h_{i_{k+1}j_{k+1}}|}{c_n}<\eps)=0\ee and suppose that Hypothesis (b)  holds with \pro tending to one. 
Then the conclusions of the theorem remain true for the convergence in probability.
\en{cor}
\bpr Recall that a sequence of real random variables converges in \pro to a deterministic limit if and only if from each of its subsequences, one can extract a subsequence converging almost surely. Hence it suffices to notice that the deterministic  theorem also holds (obviously) if one replaces the sequence $H_n$ of $n\ti n$ matrices, by a sequence $H_{\vfi(n)}$ of $\vfi(n)\ti \vfi(n)$ matrices with $\vfi(n)\to+\infty$.  
\epr

By Proposition \re{pr_crucial}, one directly deduces the following corollary.
\beg{cor}\la{tr2131212h13}Suppose that one can write $H_n=\hat{H}_n+(H_n-\hat{H}_n)$, where $H_n-\hat{H}_n$ satisfies  the hypotheses of Corollary \re{tr2131212h13StFlour} and that for a certain $\rho<1$, \be\la{hyptrunc21312} \f{\|\hat{H}_n\|}{c_n^\rho}\trm{ converges in \pro to zero.}\ee  Then the conclusions of Theorem \re{147093h48} for $H_n$ remain true for the convergence in probability.\en{cor}

\emph{Proof of Theorem \re{147093h48}.} We suppress the dependence on $n$ to simplify notation.

\noindent{\it \underline{Fact 1}:} We have $\|{H}\|_{\ell^\infty\to\ell^\infty}=|h_{i_1j_1}|(1+o(1))$ (and as a consequence, $\lam_1\le |h_{i_1j_1}|(1+o(1))$).

Indeed,    $ \|{H}\|_{\ell^\infty\to\ell^\infty}=\max_i \sum_j |h_{ij}|\ge |h_{i_1j_1}|$, thus to prove Fact 1, it suffices to notice that  for $n$ large enough,  for all $i$, $$\sum_j |h_{ij}|\le \sum_{j\ste |h_{ij}|<c_n^{\ka}}|h_{ij}|+\max_j|h_{ij}|\le c_n^{\nu}+\max_j|h_{ij}|\le \underbrace{c_n^{\nu}}_{\ll |h_{i_1j_1}|\trm{ by (a.ii).}}+|h_{i_1j_1}|.$$

\noindent{\it \underline{Fact 2}:}
For any fixed $k\ge 1$, $\lam_k\le |h_{i_kj_k}|(1+o(1))$.
 
  Indeed, the hypotheses  ensure that for $n$ large enough, the numbers $i_1, \ld, i_k$ are pairwise distinct, so that  the largest entry,  in absolute value, of the $(n-k+1)\ti (n-k+1)$ matrix 
  $H^{(k)}$,  deduced from $H$ by removing rows and columns with indices $i_1, \ld, i_{k-1}$,     is $|h_{i_kj_k}|$. This  matrix (more specifically : this \emph{sequence of matrices}, because $n$ is an implicit parameter here)  also satisfies the previous hypotheses (for the sequence $c^{(k)}_n:=c_{n+k-1}$).  Hence by the previous fact, $\lam_1(H^{(k)})\le |h_{i_kj_k}|(1+o(1))$. But by Weyl's interlacing inequalities, we have $\lam_k(H)\le \lam_1(H^{(k)})$.  It allows to conclude.

\noindent{\it \underline{Fact 3}:} For any fixed $k\ge 1$, for $\bv=\f{e^{i\tta_k}\bfe_{i_k}+e^{-i\tta_k} \bfe_{j_k}}{\sqrt{2}} $, we have $$H \bv=|h_{i_kj_k}|\bv+\br,\qquad \trm{ with  $\|\br\|=o(c_n)$}.$$

Indeed, for $\br:=H\bv-|h_{i_kj_k}|\bv$,  it is easy to see that  
$$\|\br\|\le \ff{\sqrt{2}} \lf(|h_{i_ki_k}|+| h_{j_kj_k}|+\sum_{i\notin\{ i_k, j_k\}} (|h_{i_ki}|+|h_{ij_k}|)\ri).$$
We have $\|\br\|=o(c_n)$ because by Hypothesis (b),  $|h_{i_ki_k}|+| h_{j_kj_k}|\le 2c_n^{\ka}$   and  $$\sum_{i\notin\{ i_k, j_k\}} (|h_{i_ki}|+|h_{ij_k}|)\le 2c_n^{\nu}.$$

Let us now conclude the proof of the theorem. Since $ |h_{i_kj_k}|$ has order $c_n$, Fact 3 and Part a) of Proposition \re{pr_crucial} imply that for any fixed $k\ge 1$, $H$ has an eigenvalue equal to $ |h_{i_kj_k}|(1+o(1))$. Hence by Fact 2 and Hypothesis (a.ii), $\lam_k= |h_{i_kj_k}|(1+o(1))$. By  Hypothesis (a.ii) again, it follows that $|\lam_k-\lam_{k+1}|$ has order $c_n$ and so one can apply Part b) of  Proposition \re{pr_crucial}  to deduce from fact 3 that $$\lf\|\bv_k-\f{e^{i\tta_k}\bfe_{i_k}+e^{-i\tta_k} \bfe_{j_k}}{\sqrt{2}}\ri\|\lto 0.$$\hfill$\square$

\subsection{Sums of  truncated heavy-tailed random variables}
In this section, we  give exponential estimates on the concentration of sums of truncated heavy-tailed variables. In the paper, these estimates are needed for example to give  upper bounds on the spectral radius of matrices via  the maximum of the sums of the entries along the rows.

 Let us  consider some i.i.d.   variables $Y_i\ge 0$ \st for a certain $\al>0$,    \be\la{1271221h43}\p(Y_1>y)\sL y^{-\al}\qquad\trm{as $y\to\infty$.}\ee  
 Let us also fix a sequence $d_n\sL n^\mu$ for a fixed $\mu>0$.
\beg{propo}\la{proposumHTV}
\bgt\ite[a)] For any sequence $\bet_n\sL n^b$ with $0\le b\le\mu/\al$ and any $\eps>0$, we have w.e.h.p.
$$\sum_{j=1}^{d_n}Y_j\one_{Y_j\le \bet_n}\;\le\; n^{\mu+b(1-\al)^++\eps}  .$$  
\ite[b)] If $\al\ge 1$, for any sequences $\al_n\sL n^a$,  $\bet_n\sL n^b$ with $0\le a< b\le\mu/\al$ and any $\eps>0$, we have w.e.h.p. $$\sum_{j=1}^{d_n}Y_j\one_{\al_n<Y_j\le \bet_n}\;\le\; n^{\mu-a(\al-1)+\eps}  .$$ 
\ite[c)] For any sequences $\al_n\sL n^{\f{\mu}{\al}-\eta}$, $\bet_n\sL n^{\f{\mu}{\al}+\eta'}$, with $\eta, \eta' \ge 0$ and any $\eps>\al\eta+\eta'$, we have w.e.h.p.
$$\sum_{j=1}^{d_n}Y_j\one_{\al_n<Y_j\le \bet_n}\;\le\; n^{\f{\mu}{\al}+\eps}  .$$ 
\ite[d)] For any sequences $\al_n\sL n^{\f{\mu}{\al}+\eta}$ with $\eta>0$,  $\bet_n\sL n^{\bet}$, with $\beta>0$ and any $\ga>\bet$, we have  w.e.h.p.
$$\sum_{j=1}^{d_n}Y_j\one_{\al_n<Y_j\le \bet_n}\;\le\; n^{\ga}  .$$ 
\ent
\en{propo}

Before proving the proposition, we shall first state the following   concentration result for sums of Bernoulli variables, which is a direct consequence of Bennett's inequality \cite[p. 11]{DevroyeLugosi}.
\beg{lem}\la{somme_bernoulli}
For each $n\ge 1$,  let $X_1, \ld, X_{m}$ be some independent Bernoulli  variables    with  paramater $p$ ($m$, the $X_i$'s and $p$ depending   on the parameter $n$). Suppose that $mp\ge C n^\tta$ for some constants $C, \tta>0$. Then for any fixed  $\eta>0$,  we have w.e.h.p.$$\lf|\ff{m}\sum_{i=1}^{m} X_i-p\ri|\le \eta p  .
$$
\en{lem}

 {\noindent {\it Proof of Proposition \re{proposumHTV}. }} a) First, one gets rid of the $j$'s \st $Y_j\le1$ because their sum is $\le d_n$. 
Then, set $$S_n:=\sum_{j=1}^{d_n} Y_j\one_{1<Y_j\le \bet_n}$$ and $k_\eps:=\lfloor b/\eps\rfloor$.
  We have  $ \tta:=\f{\mu-\al k_\eps \eps}{2}>0$ and for $n$ large enough, $\bet_n\le   n^{(k_\eps+1)\eps}$, so  
$$S_n\le \sum_{k=0}^{k_\eps}\sum_{j=1}^{d_n}Y_j\one_{n^{ k\eps}<Y_j\le n^{(k+1)\eps}}\le \sum_{k=0}^{k_\eps}\underbrace{\sum_{j=1}^{d_n}\one_{n^{ k\eps}<Y_j\le n^{(k+1)\eps}}}_{:= Z_{k}^{(n)}} n^{(k+1)\eps }.
 $$
 For each $k$, $Z_{k}^{(n)}$ is a sum of $d_n\sL n^\mu$ independent Bernoulli variables with parameter $p_k(n)\sL n^{-\al k\eps}$. We have $d_np_k(n)\sL n^{\mu -\al k\eps}$ hence for $n$ large enough, $d_np_k(n)\ge n^\tta$ (where $\tta$ is defined  above). As a consequence, by Lemma \re{somme_bernoulli},   w.e.h.p.,  for each $k$, $$Z_{k}^{(n)}\le 2d_np_k(n).$$ This implies that    $$ S_n\le \sum_{k=0}^{k_\eps} 2d_np_k(n)   n^{(k+1)\eps}\le n^{\mu+b(1-\al)^++\eps}.$$
 
 b) Let $S_n$ be the considered sum. The proof works in the same way as the one of a), introducing    $k_\eps:=\lfloor (b-a)/\eps\rfloor$, $\tta:=\f{\mu-\al(a+k_\eps \eps)}{2}>0$ and writing $$S_n\le  \sum_{k=0}^{k_\eps} \sum_{j=1}^{d_n}\one_{\al_nn^{ k\eps}<Y_j\le \al_nn^{(k+1)\eps}} \al_nn^{(k+1)\eps}.
 $$
 
  c) This is a direct application of Lemma \re{somme_bernoulli}, since the considered sum is $\le \bet_n \,\ti\,$the sum of $d_N$ Bernoulli variables with parameter $\sL n^{-\mu+\al\eta}$.

d)  Note that if   the considered sum is  $>n^\ga$, then there are at least $n^\ga/\bet_n$ non zero terms in the sum. By the union bound, this happens with probability at most  $$d_n^{\lceil n^\ga/\bet_n\rceil} \,\ti\, \p(Y_1>\al_n)^{\lceil n^\ga/\bet_n\rceil}$$ (indeed, there are at most $d_n^{\lceil n^\ga/\bet_n\rceil} $ subsets of cardinality ${\lceil n^\ga/\bet_n\rceil}$ in $\{1, \ld, d_n\}$). Using \eqre{1271221h43}, one can easily check that this probability is exponentially small.{\hfill $\square$\\}

 \subsection{A uniform bound for truncated moments of heavy tailed variables}  
%
%

The following result is a ``uniform in $k$ and $x$'' version of a well known result that we use  in the proof of Theorem \re{Thmoments22612}.
\beg{lem}\la{lemtruncmomentsHTRV}Let $A$ be a non negative  random variable \st for all $x\ge 0$, $$\p(A\ge x)=\ell (x)x^{-\al},$$ with $\ell$ a slowly varying function and $\al>0$. Then there exists a slowly varying function $L_0$ \st for any positive integer $k$, any $x\ge 0$, \be\la{143132}\E[A^k\one_{A\le x}]\le \beg{cases}L_0(x)&\trm{ if $k\le \al$,}\\ \\
L_0(x)\f{k}{k-\al}x^{k-\al}&\trm{ if $k> \al$.}
\en{cases}\ee
\en{lem}

\bpr First, for $k>\al$,  $$\E[A^k\one_{A\le x}]=\int_0^x kt^{k-1}\p(t\le A\le x)\ud t\le  \int_0^x kt^{k-1}\ell(t)t^{-\al}\ud t\le\tL(x)\f{k}{k-\al}x^{k-\al},$$ with $\tL(x):=\sup\{\ell(t)\ste t\in [0,x]\}$. Second, for  $k=\al$, $$\E[A^k\one_{A\le x}]\le \p(A\le 1)+\tL(x) \al\log x.$$
At last, we know, by \cite{feller2}, Chap. VIII.9, Th. 2.23, that $\E[A^k]<\infty$ when $k<\al$.  Set $$L(x):=\tL(x)+\p(A\le 1)+\tL(x) \al\log x+ \max \{\E[A^k]\ste k=1, \ld, \lfl \al\rfl\}.$$ Then  \eqre{143132} is satisfied for $L$ instead of $L_0$. Hence it suffices to prove that there is a slowly varying function $L_0$ \st 
$L\le L_0$. By Karamata's representation theorem \cite[Th. 1.3.1]{Bingham-Goldie-Teugels} (and as $\ell$ is bounded on any bounded interval), one can find a constant $C$ and a measurable function $\eps$ with null limit at $+\infty$ \st for all $t\ge 0$, $$ \ell (t)\le C\exp\int_0^t\f{\eps(u)}{u}\ud u,$$ so that $$ \tL (x)\le C\exp\int_0^x\f{|\eps(u)|}{u}\ud u,$$ which is slowly varying by Karamata's theorem again. This allows to conlude.\epr

 \begin{thebibliography}{10}
  \bibitem{ABP09} A. Auffinger, G. Ben Arous, S. P\'ech\'e  \emph{Poisson convergence for the largest eigenvalues of heavy tailed random matrices}. Ann. Inst. Henri Poincar\'e Probab. Stat. 45 (2009), no. 3, 589--610.
 \bibitem{bai-silver-book} Z.D.~Bai, J.W.~Silverstein \emph{Spectral analysis of large dimensional random matrices}. Second Edition, Springer, New York, 2009.
 \bibitem{bhatiamatrixanalysis} R. Bhatia \emph{Matrix analysis}. Graduate Texts in Mathematics, 169. Springer-Verlag, New York, 1997. 
 
  \bibitem{bhatiapertmatrices} R. Bhatia \emph{Perturbation bounds for matrix eigenvalues}. Reprint of the 1987 original. Classics in Applied Mathematics, 53. Society for Industrial and Applied Mathematics (SIAM), Philadelphia, PA, 2007.

\bibitem{GBAG} G. Ben Arous, A. Guionnet,  \emph{ The spectrum of heavy tailed random matrices.}  Comm. Math. Phys. 278 (2007), no 3, 715--751. 

\bibitem{Bingham-Goldie-Teugels} N.H. Bingham, C.M. Goldie  J.L. Teugels, \emph{Regular variation}, Cambridge University Press, 1989.

\bibitem{GB}C. Bordenave, A. Guionnet \emph{Localization and delocalization of eigenvectors for heavy-tailed random matrices}. Preprint arXiv:1201.1862.
\bibitem{CB}P. Cizeau, J-P. Bouchaud, \emph{Theory of L\'evy matrices} Phys. Rev. E 50 (1994). 
\bibitem{DevroyeLugosi} L. Devroye,  G. Lugosi \emph{Combinatorial methods in density estimation}. Springer (2001).

\bibitem{EReview}L. Erd\"os \emph{ Universality of Wigner random matrices: a Survey of Recent Results} Preprint arXiv:1004.0861. 

    \bibitem{EK2011CMP}  L. Erd\"os, A. Knowles
\emph{Quantum diffusion and eigenfunction delocalization in a random band matrix model}. 
Comm. Math. Phys. 303 (2011), no. 2, 509--554. 
    
 \bibitem{ErdKnow}  L. Erd\"os, A. Knowles, \emph{ Quantum diffusion and delocalization for band matrices with general distribution.} Ann. Henri Poincara\'e 12 (2011), no. 7, 1227--1319.   
   
   \bibitem{EKYY}  L. Erd\"os, A. Knowles, H-T. Yau, J. Yin, \emph{Delocalization and Diffusion Profile for Random Band Matrices}, arXiv:1205.5669,  to appear in Comm. Math. Phys.
   
\bibitem{ESY1}L. Erd\"os, B. Schlein, H.T. Yau \emph{Local semicircle law and complete delocalization for Wigner random
matrices}, Comm. Math. Phys. 287 (2009).
\bibitem{ESY2}L. Erd\"os, B. Schlein, H.T. Yau, \emph{Semicircle law on short scales and delocalization of eigenvectors for Wigner random matrices}, Ann. Prob. 37 (2009).
\bibitem{ESY3}L. Erd\"os, B. Schlein, H.T. Yau \emph{Wegner estimate and level repulsion for Wigner random matrices}, Int. Math. Res. Not. 2010 (2010).

 \bibitem{feller2} W. Feller \emph{An introduction to probability theory and its applications}, volume II, second edition, New York London Sydney : J. Wiley, 1966.
\bibitem{FM}Y.V. Fyodorov and A.D. Mirlin, \emph{Scaling properties of localization in random band matrices: a $\sigma$-model approach}, Phys. Rev. Lett. 67 (1991), no. 18, 2405-2409.
 
\bibitem{KY}A. Knowles, J. Yin, \emph{Eigenvector Distribution of Wigner Matrices}, arXiv:1102.0057,  to appear in Prob. Theor. Rel. Fields.

\bibitem{KY2012} A. Knowles, J. Yin, \emph{The Outliers of a Deformed Wigner Matrix}, arXiv:1207.5619, to appear in Ann. Prob.

\bibitem{LLRextremes} M.R. Leadbetter, G.Lindgren and H. Rootz\'en \emph{Extremes and Related Properties of Random Sequences and Processes}, Springer-Verlag, New York, 1983.
\bibitem{soshnipeche07} S. P\'ech\'e, A. Soshnikov  \emph{Wigner random matrices with non-symmetrically distributed entries}. J. Stat. Phys. 129 (2007), no. 5-6, 857--884.
\bibitem{resnick} R. Resnick \emph{Extreme Values, Regular Variation and Point Processes}. Springer, New York, 1987.
\bibitem{Schenker} J. Schenker 
\emph{Eigenvector localization for random band matrices with power law band width. }
Comm. Math. Phys. 290 (2009), no. 3, 1065–1097. 
\bibitem{soshnisinai98} Y. Sinai, A. Soshnikov  \emph{A refinement of Wigner's semicircle law in a neighborhood of the spectrum edge for random symmetric matrices}. (Russian) Funktsional. Anal. i Prilozhen. 32 (1998), no. 2, 56--79, 96; translation in Funct. Anal. Appl. 32 (1998), no. 2, 114--131.
\bibitem{Sodin} S. Sodin \emph{The spectral edge of some random band matrices}, Ann. Math. 172 (2010), 2223--2251. 
\bibitem{Sos}Y. Sinai, A. Soshnikov \emph{ Central limit theorem for traces of large random symmetric matrices}. Bol.
Soc. Brasil. Mat. 29 (1998) No. 1, 1–24.
\bibitem{soshni99} A. Soshnikov  \emph{Universality at the edge of the spectrum in Wigner random matrices}. Comm. Math. Phys. 207 (1999), no. 3, 697--733.
\bibitem{sosh04}  A. Soshnikov  \emph{Poisson statistics for the largest eigenvalues of Wigner random matrices with heavy tails}. Electron. Comm. Probab. 9 (2004), 82--91.
\bibitem{sosh06}  A. Soshnikov  \emph{Poisson statistics for the largest eigenvalues in random matrix ensembles}. Mathematical physics of quantum mechanics, 351--364, Lecture Notes in Phys., 690, Springer, Berlin, 2006.

\bibitem{Spencer} T. Spencer \emph{Random banded and sparse matrices} The Oxford handbook of random matrix theory. Oxford University Press (2011) 471--488.
\bibitem{TV1}T. Tao, V. Vu \emph{Random matrices: Universal properties of eigenvectors}, Random Matrices: Theory Appl. 1 (2012) 1150001.
\bibitem{TV2}T. Tao, V. Vu \emph{Random matrices: universality of local eigenvalue statistics up to the edge}. 
Comm. Math. Phys. 298 (2010), no. 2, 549--572.
 
 \en{thebibliography}
\en{document}